\documentclass[a4paper,12pt]{amsart}
\usepackage{mathrsfs}
\usepackage{amsfonts}
\usepackage{amssymb}
\usepackage{amscd}
\usepackage{ifthen}
\usepackage{graphicx}
\usepackage[usenames]{color}

\nonstopmode \numberwithin{equation}{section}
\newtheorem{thm}{Theorem}
\newtheorem{cor}{Corollary}
\newtheorem{lem}{Lemma}

\newtheorem{claim}{Claim}
\newtheorem{conj}[equation]{Conjecture}

\theoremstyle{definition}
\newtheorem{defn}{Definition}
\newtheorem{case}{Case}
\newtheorem{subcase}{Subcase}

\newtheorem{examp}{Example}
\newtheorem{prob}[equation]{Problem}
\newtheorem{ques}[equation]{Question}
\newtheorem{rem}{Remark}

\newcounter {own}
\def\theown {\thesection       .\arabic{own}}

\newenvironment{pf}[1][]{%
 \vskip 3mm
 \noindent
 \ifthenelse{\equal{#1}{}}%
  {{\slshape Proof. }}%
  {{\slshape #1.} }%
 }%
{\qed\bigskip}

\newcounter{alphabet}
\newcounter{tmp}
\newenvironment{Thm}[1][]{\refstepcounter{alphabet}%
\bigskip%
\noindent%
{\bf Theorem \Alph{alphabet}}%
\ifthenelse{\equal{#1}{}}{}{ (#1)}%
{\bf .} \itshape}{\vskip 8pt}

\makeatletter
\newcommand{\Ref}[1]{\@ifundefined{r@#1}{}{\setcounter{tmp}{\ref{#1}}\Alph{tmp}}}
\makeatother

\newenvironment{Lem}[1][]{\refstepcounter{alphabet}%
\bigskip%
\noindent%
{\bf Lemma \Alph{alphabet}}%
{\bf .} \itshape}{\vskip 8pt}

\newcommand{\IC}{{\mathbb C}}
\newcommand{\ID}{{\mathbb D}}

\def\be{\begin{equation}}
\def\ee{\end{equation}}

\newcommand{\bee}{\begin{enumerate}}
\newcommand{\eee}{\end{enumerate}}

\newcommand{\blem}{\begin{lem}}
\newcommand{\elem}{\end{lem}}
\newcommand{\bthm}{\begin{thm}}
\newcommand{\ethm}{\end{thm}}
\newcommand{\bcor}{\begin{cor}}
\newcommand{\ecor}{\end{cor}}
\newcommand{\beg}{\begin{examp}}
\newcommand{\eeg}{\end{examp}}
\newcommand{\begs}{\begin{examples}}
\newcommand{\eegs}{\end{examples}}
\newcommand{\bdefe}{\begin{defn}}
\newcommand{\edefe}{\end{defn}}
\newcommand{\bprob}{\begin{prob}}
\newcommand{\eprob}{\end{prob}}
\newcommand{\bques}{\begin{ques}}
\newcommand{\eques}{\end{ques}}
\newcommand{\bei}{\begin{itemize}}
\newcommand{\eei}{\end{itemize}}
\newcommand{\bca}{\begin{case}}
\newcommand{\eca}{\end{case}}
\newcommand{\bsca}{\begin{subcase}}
\newcommand{\esca}{\end{subcase}}

\newcommand{\bcl}{\begin{claim}}
\newcommand{\ecl}{\end{claim}}

\newcommand{\bcon}{\begin{conj}}
\newcommand{\econ}{\end{conj}}
\newcommand{\bcons}{\begin{conjs}}
\newcommand{\econs}{\end{conjs}}
\newcommand{\bprop}{\begin{propo}}
\newcommand{\eprop}{\end{propo}}
\newcommand{\br}{\begin{rem}}
\newcommand{\er}{\end{rem}}
\newcommand{\brs}{\begin{rems}}
\newcommand{\ers}{\end{rems}}
\newcommand{\bo}{\begin{obser}}
\newcommand{\eo}{\end{obser}}
\newcommand{\bos}{\begin{obsers}}
\newcommand{\eos}{\end{obsers}}
\newcommand{\bpf}{\begin{pf}}
\newcommand{\epf}{\end{pf}}
\newcommand{\ba}{\begin{array}}
\newcommand{\ea}{\end{array}}
\newcommand{\beq}{\begin{eqnarray}}
\newcommand{\beqq}{\begin{eqnarray*}}
\newcommand{\eeq}{\end{eqnarray}}
\newcommand{\eeqq}{\end{eqnarray*}}

\newcommand{\ds}{\displaystyle}

\newcounter{minutes}
\divide\time by 60
\newcounter{hours}
\multiply\time by 60 \addtocounter{minutes}{-\time}

\begin{document}
\bibliographystyle{amsplain}
\title {Extreme points and support points of families of harmonic Bloch mappings}

\thanks{
File:~\jobname .tex,
          printed: \number\day-\number\month-\number\year,
          \thehours.\ifnum\theminutes<10{0}\fi\theminutes}

\author[H. Deng]{Hua Deng}
\address{H. Deng, Department of Mathematics, Hebei University,
Baoding, Hebei 071002, People's Republic of China.}
\email{1120087434@qq.com}

\author[S. Ponnusamy]{Saminathan Ponnusamy}

\address{S. Ponnusamy, Department of Mathematics, Indian Institute of
Technology Madras, Chennai-600 036, India. }
\email{samy@iitm.ac.in}

\author{Jinjing Qiao${}^{~\mathbf{*}}$}
\address{J. Qiao, Department of Mathematics, Hebei University,
Baoding, Hebei 071002, People's Republic of China.}
\email{mathqiao@126.com}

\keywords{Bloch function, Harmonic Bloch mapping,  extreme point, support point.\\
${}^{\mathbf{*}}$ Corresponding author
}


\begin{abstract}
In this paper, the main aim is to discuss the existence of the
extreme points and support points of families of harmonic Bloch mappings and
little harmonic Bloch mappings. First, in terms of the Bloch unit-valued
set, we prove a necessary condition for a  harmonic Bloch mapping
(resp. a little harmonic Bloch mapping) to be an extreme point of
the unit ball of the  normalized harmonic Bloch
 spaces (resp. the normalized little harmonic Bloch
 spaces) in the unit disk $\mathbb{D}$. Then we show
that a  harmonic  Bloch mapping $f$ is a support point of the unit
ball of the normalized harmonic Bloch spaces in $\mathbb{D}$ if and
only if the Bloch unit-valued set of $f$ is not empty. We also give
a characterization for the support points of the unit ball of the
harmonic Bloch spaces in $\mathbb{D}$.
\end{abstract}

\thanks{}

\maketitle \pagestyle{myheadings} \markboth{H. Deng, S. Ponnusamy
and J. Qiao}{Extreme points and support points of families of harmonic Bloch
mappings}

\section{Introduction and Preliminaries }\label{csw-sec1}
Support points and extreme points of analytic functions play important roles
in solving extremal problems. It is known that in the topology of uniform
convergence on compacta, any compact family of analytic
functions contains support points and the set of all support
points contains an extreme point. This remarkable fact plays an active role in
solving extremal problems for various families of analytic functions
(see \cite{Am,Am1,BCK,BHHW,Du2,HH,HM,Mac,Mac1,RW,ST,Tk,WirthXiao} and for very recent work on this
topic, we refer to \cite{KW-2019,KW-2019b}).
The main focus in this article is to extend a number of results
from the theory of  analytic functions to the case of planar harmonic mappings.
In particular, we extend the work of Cima and Wogen \cite[Theorem 2]{CW} in
the setting of little harmonic Bloch mappings, and construct a counterexample to show that
\cite[Theorem 1 and Corollary 1]{CW} fail to hold for (little) harmonic Bloch mappings.
Moreover, we establish a characterization for a harmonic Bloch mapping to be a support point
of $\mathscr{B}_{\mathcal {H}, 1}$ which in turn extends the work of Bonk \cite[Theorem 3]{Bo}.
The definitions of these mappings and the exact formulation of the results of Cima and Wogen
will be addressed later in this section and the results of Bonk in the next section.

Let  $\IC$ be the complex plane, and  $\Omega$ be a simply connected domain in $\IC$. A
harmonic mapping $f$ on $\Omega$ is a complex-valued function of the form $f=u+iv$, where $u$ and $v$ are
real-valued harmonic functions on $\Omega$. This function has the canonical decomposition $f=h+\overline{g}$,
where $h$ and $g$ are analytic functions in $\Omega $, known as analytic and co-analytic parts of $f$, respectively,
and $g(z_0)=0$ for some prescribed point $z_0\in \Omega$.

In the following, we introduce some necessary notions and notations. Let $\mathbb{D}_r=\{z\in \mathbb{C}:|z|<r\}$
for $r>0$. Throughout this paper, we
consider harmonic mappings in the unit disk $\mathbb{D}=\mathbb{D}_1$.
Let $\mathcal {H}(\mathbb{D})$ denote the class of all harmonic mappings
in $\mathbb{D}$ and $\mathcal {A}(\mathbb{D})$ the set of all analytic functions in $\mathbb{D}$.

A function $f=h+\overline{g}\in \mathcal {H}(\mathbb{D})$ is called a {\it harmonic Bloch mapping} if
$$ \beta_f:=\sup_{z\neq w}\frac{|f(z)-f(w)|}{\rho(z,w)}<\infty,
$$
where $\beta_f$ is called the ${\it Bloch\; constant}$ of $f$ and
$$\rho(z,w)=\frac{1}{2}\log\left(\frac{1+|\frac{z-w}{1-\overline{z}w}|}
{1-|\frac{z-w}{1-\overline{z}w}|}\right)={\rm arctanh}
\Big|\frac{z-w}{1-\overline{z}w}\Big|
$$
denotes the hyperbolic distance between $z$ and $w$ in $\mathbb{D}$ (cf. \cite{Co}).
Moreover, it is known that
$$\beta_{f}=\sup_{z\in\mathbb{D}}\mu_f(z), \quad
\mu_f(z):=(1-|z|^2)\big(|h'(z)|+|g'(z)|\big).
$$
Obviously, the correspondence $f\mapsto \beta_f$ is invariant under pre-composition by conformal automorphisms of $\mathbb{D}$. We remark that in the
case of an analytic function $f$, its Bloch constant (see \cite{ACP} and \cite[Theorem 10]{Co1} for details)
naturally takes the form
$$\beta_f=\sup_{z\in \mathbb{D}}(1-|z|^2)|f'(z)|
$$
and $f$ is a Bloch function if $\beta_f<\infty$.
Furthermore, a harmonic mapping $f=h+\overline{g}$ is said to be Bloch if and only if both $h$ and
$g$ are (analytic) Bloch functions. This can be seen from the fact
$$\max\big\{\beta_h, \beta_g\big\}\leq\beta_f\leq \beta_h+ \beta_g.
$$

Let $\mathscr{B}_\mathcal {H}$ (resp. $\mathscr{B}$) denote the
class of all harmonic mappings (resp. analytic functions) $f$ with $\beta_f<\infty$. It is easy to see
that $\mathscr{B}_\mathcal {H}$ (resp. $\mathscr{B}$) is a {\it
Banach space} with the norm
$$\|f\|=|f(0)|+\beta_f,
$$
which is called the {\it harmonic  $(${\rm resp.} analytic$)$ Bloch
space}. Each element in $\mathscr{B}_\mathcal {H}$ (resp.
$\mathscr{B}$) is  a {harmonic Bloch mapping} (resp. a {Bloch
function}).

The {\it little harmonic $(${\rm resp.} analytic$)$ Bloch space}
$\mathscr{B}_{\mathcal {H},0}$ (resp. $\mathscr{B}_{0}$) is the set
of all mappings $f\in \mathscr{B}_\mathcal {H}$ (resp. $f\in
\mathscr{B}$) satisfying
$$\lim_{|z|\rightarrow 1}\mu_f(z)=0.
$$
Each element in  $\mathscr{B}_{\mathcal {H},0}$ (resp. $\mathscr{B}_{0}$) is called a {\it little
harmonic Bloch mapping} (resp. a {\it little
Bloch function}). Also we let
\beqq
\mathscr{B}_{\mathcal {H}, 1}&=&\{f\in \mathscr{B}_\mathcal {H}: \|f\|\leq 1\},\\
\mathscr{B}_{\mathcal {H},0, 1}&=&\{f\in \mathscr{B}_{\mathcal {H},0}: \|f\|\leq 1\},\\
\widetilde{\mathscr{B}}_{\mathcal {H}, 1}&=& \{f\in \mathscr{B}_{\mathcal {H}, 1}: h(0)=g(0)=0\},\\
\widetilde{\mathscr{B}}_{\mathcal {H}, 0,1} &=& \{f\in \mathscr{B}_{\mathcal {H},0,1}: h(0)=g(0)=0\}, ~\mbox{ and }\\
\Lambda_f &=&\{z\in \mathbb{D}: \mu_f(z)=1\}.
\eeqq
In particular, $\Lambda_f$ is called the {\it Bloch unit-valued set} of $f$.
It is natural to set
$$\left\{ \ba{ll}\mathscr{B}_{1}=\mathscr{B}_{\mathcal {H}, 1}\cap\mathcal
{A}(\mathbb{D}), &~\widetilde{\mathscr{B}}_{1}= \widetilde{\mathscr{B}}_{\mathcal {H}, 1}\cap\mathcal{A}(\mathbb{D})  ~\\
\mathscr{B}_{0,1}=\mathscr{B}_{\mathcal {H},0, 1}\cap\mathcal
{A}(\mathbb{D}), &~ \widetilde{\mathscr{B}}_{0,1}= \widetilde{\mathscr{B}}_{\mathcal {H}, 0, 1}\cap\mathcal
{A}(\mathbb{D}).
\ea \right .
$$

\bdefe\label{de1.1}
Let $X$ be a topological vector space over the
field of complex numbers, and let $D$ be a convex subset of $X$. A point $x\in
D$ is called an {\it extreme point} of $D$ if it has no
representation of the form $x=ty+(1-t)z$ $(0<t<1)$ as a proper
convex combination of two distinct points $y$ and $z$ in $D$. A
point $x\in D$ is called a \textit{support point} of $D$ if there is
a continuous linear functional $J$, not constant on $D$, such that
${\rm Re\,}\{J(x)\}\geq {\rm Re\,}\{J(y)\}$ for all $y\in D$ (cf. \cite{Du2}).
\edefe


For (analytic) Bloch functions, in \cite{CW}, it is shown that the set of all
extreme points of the unit ball ${\mathscr{B}}_{0,1}$ in the (analytic) little Bloch space in $\mathbb{D}$ is the union of the
set of all unimodular constants and the set of extreme points of the
convex set $\widetilde{\mathscr{B}}_{0,1}$, which is compact in the topology of uniform convergence on compacta
(see \cite[Corollary $2$]{CW}).
The authors in \cite{CW} also proved that a sufficient condition for a
function $f\in \widetilde{\mathscr{B}}_{1}$ to be an extreme point of
$\widetilde{\mathscr{B}}_{1}$ is that the intersection of
$\Lambda_f$ with the disk $\overline{\ID}_R:=\{z\in \mathbb{C}:\; |z|\leq R\}$ for some $0<R< 1$ has to be infinite
 (see \cite[Theorem $1$]{CW}); Further they obtained that, under the assumption
$\lim_{|z|\rightarrow 1}\mu_f(z)=0$, the condition ``$\Lambda_f$ being infinite" is
necessary for $f$ to be extreme (see \cite[Theorem $2$]{CW}). In \cite{Bo}, a characterization of support points in
${\mathscr{B}}_{1}$ was established  in terms of the set $\Lambda_f$, see \cite[Theorem $3$]{Bo}.

The main aim of this paper is to extend the results stated as above to the
case of harmonic mappings and the results are organized as follow.  In Section \ref{sec-2}
(see Theorem \ref{thm2.0}), we prove that \cite[Theorem $2$]{CW} holds for
the setting of little harmonic Bloch mappings.  Then we construct a counterexample to show that
\cite[Theorem $1$ and Corollary $1$]{CW} fail to hold for (little) harmonic Bloch mappings.
In the end, we consider the support points of ${\mathscr{B}}_{\mathcal {H}, 1}$. In Section \ref{sec-3}
(see Theorem \ref{thm2.1}),   a characterization for a harmonic Bloch mapping to be a support point
of $\mathscr{B}_{\mathcal {H}, 1}$ is established, and this result is indeed a generalization of
\cite[Theorem 3]{Bo}.

\section{Extreme  points }\label{sec-2}

\subsection{The statement of the first main result}


\bthm\label{thm2.0}
\bee
\item[{\rm (1)}]  Suppose that $f\in \widetilde{\mathscr{B}}_{\mathcal {H},0,1}$ and that $f$ is an extreme point
of $\widetilde{\mathscr{B}}_{\mathcal {H},0,1}$. Then $\Lambda_f$ is
infinite.

\item[{\rm (2)}] Suppose that $f\in \widetilde{\mathscr{B}}_{\mathcal {H},1}$ and that
$f$ is an extreme point of $\widetilde{\mathscr{B}}_{\mathcal {H},1}$.
Then there exists an $R\in (0,1)$ such that the intersection $\Lambda_f\cap\{z:
|z|<R\}$ is an  infinite set.
\eee
\ethm

We remark that if $f=h+\overline{g}\in \widetilde{\mathscr{B}}_{\mathcal
{H},0,1}$ and $\Lambda_h$ is infinite, then $f=h$ and so \cite[Corollary $1$]{CW} implies that $f$ is an
extreme point of $\widetilde{\mathscr{B}}_{\mathcal {H},0,1}$.

\subsection{Example}
We demonstrate by an example that Theorem~1 and Corollary~1 \cite{CW} fail to hold for the
corresponding class of harmonic Bloch mappings.



For $a\in (0,2)$, consider $f_a(z)=h_a(z)+\overline{g_{a}(z)}$, where
$$ h_a(z)=\frac{3\sqrt{3}}{8}az^2 ~\mbox{ and }~ g_a(z)=-h_{2-a}(z).
$$
Then for each $a\in (0,2)$, we have
$$|h_a'(z)|+|g_a'(z)| =\frac{3\sqrt{3}}{2}|z| ~\mbox{ and }~\mu_{f_a}(z)=\frac{3\sqrt{3}}{2}|z|(1-|z|^2).
$$
Since
$$\sup_{|z|<1} \mu_{f_a}(z) \, = \, \frac{3\sqrt{3}}{2} \sup_{0\leq x<1} x(1-x^2) = 1,
$$
it follows that $f_a\in \widetilde{\mathscr{B}}_{\mathcal {H},0,1}$ for each $a\in (0,2)$.
Also, $ f_a\in \widetilde{\mathscr{B}}_{\mathcal {H},1}$ for each $a\in (0,2)$. Moreover, $\mu_{f_a}(z)=1$ if and only if $A(|z|)=0$, where
$$A(r)=r^3-r +\frac{2}{3\sqrt{3}}, \quad r\in [0,1).
$$
As $A'(r)=3r^2-1$, the only critical point on $(0,1)$ is at $r=1/\sqrt{3}$ and it follows easily that
$A(r)$ is decreasing on $[0,1/\sqrt{3})$ and increasing on $(1/\sqrt{3},1]$. Consequently, $A(r)\geq A(1/\sqrt{3})=0$ for all
$r\in [0,1)$ and thus, $\Lambda_{f_a}=\{z\in \ID:\, |z|=1/\sqrt{3}\}$. Hence,  $\Lambda_{f_a}$ is infinite.
Finally, it is a simple exercise to see that for  $a\in (0,2)\backslash\{1\}$,
$$f(z)=\frac{1}{2}\big(f_{a}(z)+f_{2-a}(z)\big)
$$
which implies that  $f$ is neither an extreme point of $\widetilde{\mathscr{B}}_{\mathcal {H},0,1}$ nor that of
$\widetilde{\mathscr{B}}_{\mathcal {H},1}$.


%
%
%

\subsection{Several lemmas}

The proof of Theorem \ref{thm2.0} is based on several lemmas.
Before the statement of our first lemma, let us record a result from \cite[p.~145]{ZS}
(see also \cite[Theorem A]{CW}) which is a real analytic version of the Weierstrass Preparation
Theorem.

\begin{Thm}\label{Lem2.0}
Let $G(x,y)$ be a convergent real power series such that $G(0,0)=0$ and
$G(0,y)=\sum_{n=s}^{\infty}b_ny^n$, where $s\geq 1$ and $b_s \neq
0$. Then there are power series $\Omega(x,y)$, $A_i(x)$ $(i=0,
\ldots, s-1)$ such that
$$ G(x,y)=(y^s+ A_{s-1}(x)y^{s-1}+\cdots+A_0(x))\Omega(x,y)
$$
and $\Omega(0,0)\neq 0$.
\end{Thm}

Let us now state and prove our first lemma.

\blem\label{lem3.2}
Suppose that $f=h+\overline{g}\in \widetilde{\mathscr{B}}_{\mathcal {H},1}$ with $|h'(0)|=1$ or
$|g'(0)|=1$, and that there is a $\delta_0>0$ satisfying
\be\label{eq-ass1}
\big(|h'(z)|+|g'(z)|\big)(1-|z|^2)<1 ~\mbox{ for $0<|z|<\delta_0$}.
\ee
Then there exists a positive integer $n$, and a $\delta\in (0,\delta_0]$ such that
\be\label{eq-ass1a}
\big(|h'(z)|+|g'(z)|+|z|^n\big)(1-|z|^2)<1  ~\mbox{ for $0<|z|<\delta$}.
\ee
\elem

\bpf  Without loss of generality, we may assume that $|h'(0)|=1$. By considering the function
$e^{i\theta_1}h+e^{i\theta_2}\overline{g}$, if needed, we assume further
that $h'(0)=1$. Then $h'$ and $g'$ have the following series expansions:
$$
h'(z)=1+a_1z+a_2z^2+\sum_{k=3}^{\infty}a_kz^k ~\mbox{ and }~
g'(z)=b_0+b_1z+b_2z^2+\sum_{k=3}^{\infty}b_kz^k.
$$
Since $h$ and $g$ are (analytic) Bloch functions, it follows from
the similar reasoning as in the proof of \cite[Lemma 2]{CW} that
$a_1=b_0=b_1=0$, $|a_2|\leq 1$ and $|b_2|\leq 1$, so that
\be\label{eq-ass2}
h'(z)=1+a_2z^2+\sum_{k=3}^{\infty}a_kz^k ~\mbox{ and }~ g'(z)=b_2z^2+\sum_{k=3}^{\infty}b_kz^k.
\ee

\bcl\label{cl-1} $|a_2|+|b_2|\leq 1$.
\ecl

Suppose on the contrary that $|a_2|+|b_2|>1$. Then both $a_2$ and $b_2$ must be non-zero.
With $z=re^{i\theta}$, we can choose a suitable $\theta$ so that for all
sufficiently small $r$ we have $|1+a_2z^2|=1+|a_2|r^2$ and
$$\left |\sum_{k=3}^{\infty}a_kz^k\right |\leq \left (\frac{|a_2|+|b_2|-1}{4}\right )r^2,
$$
because  $\sum_{k=3}^{\infty}a_kz^k$ is continuous at the origin. Similarly we have
$$ \left |b_2z^2+\sum_{k=3}^{\infty}b_kz^{k}\right |=|b_2z^2|
\left |1+\sum_{k=3}^{\infty}\frac{b_k}{b_2}z^{k-2}\right |\geq
|b_2z^2|\left (1-\frac{|a_2|+|b_2|-1}{4|b_2|}\right )
$$
and therefore, using these two inequalities, we deduce that
\beq\label{eq2.17}
|h'(re^{i\theta})|+|g'(re^{i\theta})|
&=&\nonumber \left |1+a_2z^2+\sum_{k=3}^{\infty}a_kz^k\right |+ \left|b_2z^2+\sum_{k=3}^{\infty}b_kz^k \right |\\
&\geq& \nonumber\left |1+a_2z^2\right |-\left |\sum_{k=3}^{\infty}a_kz^k \right|+|b_2z^2| \left (1-\frac{|a_2|+|b_2|-1}{4|b_2|}\right )\\
 &\geq&\nonumber
1+ \left(|a_2|-\frac{|a_2|+|b_2|-1}{4}\right )r^2 +\left(|b_2|-\frac{|a_2|+|b_2|-1}{4} \right )r^2\\
&=&  1+\left (\frac{|a_2|+|b_2|+1}{2}\right )r^2.
\eeq
Since $h'(0)=1$ and $b_0=g'(0)=0$, we know that $\|f\|=1$, which contradicts with \eqref{eq2.17}.
Hence Claim \ref{cl-1} holds.

\medskip

Based on Claim \ref{cl-1}, we divide the rest of the proof into two cases.

\bca $|a_2|+|b_2|<1 $.\eca

In this case it is obvious from the continuity that, there is a
$\delta_1\in (0, \delta_0]$ so that for all $z$ with
$0<|z|<\delta_1$, we have \be\label{eq-Az} A(z)= |z|^3+
\left|\sum_{k=3}^{\infty}a_kz^k \right|+
\left|\sum_{k=3}^{\infty}b_kz^k\right|<(1-|a_2|-|b_2|)|z|^2 \ee and
thus, by \eqref{eq-Az}, we obtain that \beq\label{qw-1}
|h'(z)|+|g'(z)|+|z|^3 &\leq &1+|a_2|\,|z|^2+ |b_2|\,|z|^2 +A(z)\\
\nonumber &<& 1+|z|^2,
\eeq
which shows that
$$ (|h'(z)|+|g'(z)|+|z|^3)(1-|z|^2)<1-|z|^4<1 ~\mbox{ for $0<|z|<\delta_1$}.
$$
Thus, in this case, \eqref{eq-ass1a} holds with $n=3$ and $\delta =\delta_1$.

\bca $|a_2|+|b_2|=1 $.\eca

In this case, we need to deal with three subcases separately.

\bsca $0<|b_2|<1$.\esca

Obviously, $a_2\not=0$. By using a rotation, we may assume that $a_2> 0$. Now, we let
$$G(x,y)=\big(1-(x^2+y^2)\big)^{-2}-\big(|h'(x+iy)|+|g'(x+iy)|\big)^2,
$$
where $z=x+iy$. Clearly, $G(0,0)=0$ and, by the assumption \eqref{eq-ass1}, we have
\be\label{eq2.18}
G(x,y)>0~\mbox{ for $0<|z|<\delta_0$}.
\ee
As with standard practice, we denote a convergent power series having only terms of order $n$
$(n\geq 1)$ or higher by $O_n$. Accordingly,
$$|h'(x+iy)|^2=1+2a_2(x^2-y^2)+O_3(x,y)~\mbox{ and }~ |g'(x+iy)|^2=O_4(x,y).
$$
It follows from $0<|b_2|<1$ and \eqref{eq-ass2} that
$$2|h'(x+iy)|\,|g'(x+iy)|=2|b_2  z^2|\, |\varphi(z)|,
$$
where
$$\varphi(z)=\sum_{k=2}^{\infty}\frac{b_k}{b_2}z^{k-2}+
\left (\sum_{k=2}^{\infty}\frac{b_k}{b_2}z^{k-2}\right )\left (\sum_{k=2}^{\infty}a_kz^k\right ).
$$
Note that $\varphi$ is analytic in $|z|<\delta_2'$ for a  $\delta_2'\in
(0, \delta_0]$. Since $\varphi(0)=1\neq0$, there exists a $\delta_2''\in (0, \delta_2']$ such that $\varphi(z)\neq0$ for
$|z|<\delta_2''$.  Since a non-vanishing analytic function in $|z|<\delta_2''$ admits a square root, there exists
an analytic function $\psi$ such that $\psi(z)^2=\varphi(z)$ and $\psi$ has the expression of the form
$$\psi(z)=1+\sum_{k=1}^{\infty}c_kz^k,
$$
and so,
$|\varphi(z)|=|\psi(x+iy)|^2=1+O_1(x,y).
$
Hence
$$ 2|h'(x+iy)|\,|g'(x+iy)|=2|b_2|\,|z|^2+O_3(x,y),
$$
whence
$$G(x,y)=2(1+a_2-|b_2|)y^2+O_3(x,y).
$$
Therefore the similar reasoning as in the proof of \cite[Lemma 1]{CW} and \eqref{eq2.18} shows that there is
$\delta_2'''$,  $0<\delta_2'''\leq \delta_2''$ and a positive integer $n_0$ such
that for each pair $x$ and $y$ with $0<|z|<\delta_2'''$,
\be\label{eq2.2}
|z|^{2n_0}=(x^2+y^2)^{n_0}<G(x,y).
\ee
Let $n_1=2n_0+1$,
$$K_0=\sup_{f\in{\mathscr{B}}_{\mathcal
{H},1}}\{2|h'(z)|+2|g'(z)|+|z|^{2n_0}: \;|z|<\delta_0\},  ~\mbox{ and }~\delta_2=\min\Big\{\delta_2''',\frac{1}{{K_0}}\Big\}.
$$
Then, by \eqref{eq2.2}, we see that for $z$ with $0<|z|<\delta_2$,
\beq\label{qw-2}
(|h'(z)|+|g'(z)|+|z|^{n_1})^2
&=&(|h'(z)|+|g'(z)|)^2+|z|^{2n_1}\\ \nonumber
&& + \,2|z|^{n_1}(|h'(z)|+|g'(z)|)\\ \nonumber
&\leq&(|h'(z)|+|g'(z)|)^2+K_0|z|^{n_1}\\ \nonumber
&<&\frac{1}{(1-|z|^2)^2}
\eeq
and thus, \eqref{eq-ass1a} holds with $n=n_1$ and $\delta =\delta_2$.

\bsca $b_2=0$.\esca

It follows that $|a_2|=1$. Without loss of generality, we may assume $a_2=1$. Again there are two cases.

At first, if $g'(z)=\sum_{k=2m}^{\infty}b_kz^k$, where $m> 1$ and $b_{2m}\neq 0$,
then the similar reasoning as in the discussions of Subcase $1$
shows that there exist a $\delta_3'>0$ and an $n_2$ such that
\beq\label{qw-3}
\big(|h'(z)|+|g'(z)|+|z|^{n_2}\big)(1-|z|^2)<1 ~\mbox{ for $0<|z|<\delta_3'$}.
\eeq
and thus, \eqref{eq-ass1a} holds with $n=n_2$ and $\delta =\delta_3'$.

Secondly, if $g'(z)=\sum_{k=2m+1}^{\infty}b_kz^k$, where $m\geq 1$ and
$b_{2m+1}\neq 0$, the assumption \eqref{eq-ass1} and \eqref{eq-ass2} tell
us  that for $z$ with $0<|z|<\delta_0$,
\be\label{eq2.3}
\left |1+\sum_{k=2}^{\infty}a_kz^k
\right |+|z|^{2m+1}\,\left |b_{2m+1}+\sum_{k=2m+2}^{\infty}b_{k}z^{k-2m-1}
\right|< 1+|z|^2+\sum_{k=2}^{\infty}|z|^{2k}.
\ee
Obviously, there exists a $\delta_3''$, $0<\delta_3''<\delta_0$, such that
$$\left |1+\sum_{k=2}^{\infty}a_kz^k\right |\neq
0 ~\mbox{ and }~ \left |b_{2m+1}+\sum_{k=2m+2}^{\infty}b_{k}z^{k-2m-1}\right |\neq0 ~\mbox{ for $0<|z|<\delta_3''$}.
$$
It follows from \eqref{eq2.3} that for $z$ with $0<|z|<\delta_3''$,
\be\label{eq2.4}
|z|^{2m+1}\left |\frac{b_{2m+1}+\sum_{k=2m+2}^{\infty}b_{k}z^{k-2m-1}}{1+\sum_{k=2}^{\infty}a_kz^k}\right |
<\frac{1+|z|^2+\sum_{k=2}^{\infty}|z|^{2k}}{|1+\sum_{k=2}^{\infty}a_kz^k|}-1.
\ee

Obviously, $\Phi (z)=\frac{1}{1+\sum_{k=2}^{\infty}a_kz^k}$ is non-vanishing and analytic in
$|z|<\delta_3''$ and therefore, there exists an analytic function $\Psi$ such that $\Phi(z)=\Psi ^2(z)$ in
$|z|<\delta_3''$. As $\Phi(0)=1$, $\Phi'(0)=0$ and $\Phi'' (0)/2\, =-a_2=-1$, we find that
$\Psi (0)=1$,  $\Psi'(0)=0$ and $\Psi'' (0)/2\, =-a_2/2\,=-1/2$. Thus, we have the following series expansion for $\Phi$ and
$\Psi$:
$$\Phi(z)=  1-z^2 +a_3^{*}z^3 +\cdots =  \Psi ^2(z),
\quad \Psi (z)=1-\frac{z^2}{2} +b_3^{*}z^3 +\cdots.
$$
Clearly, the last relation yields
$$|\Phi(z)|=\Psi (z)\overline{\Psi(z)}=\left (1-\frac{z^2}{2} +b_3^{*}z^3 +\cdots\right )\overline{\left ( 1-\frac{z^2}{2} +b_3^{*}z^3 +\cdots\right)}
$$
from which we obtain that
$$\frac{1}{|1+\sum_{k=2}^{\infty}a_kz^k|}=1-{\rm Re}(z^2)+O_3(x,y)=1-(x^2-y^2)+O_3(x,y).
$$
It follows that
\be\label{eq-Rxy}
R(x,y):=\frac{1+|z|^2+\sum_{k=2}^{\infty}|z|^{2k}}{|1+\sum_{k=2}^{\infty}a_kz^k|}-1=2y^2+O_3(x,y).
\ee
By Theorem \Ref{Lem2.0}, we have
\be\label{eq-Rxy-extra}
R(x,y)=\big(y^2+A_1(x)y+A_0(x)\big)\Omega(x,y),
\ee
where the functions $A_0$, $A_1$ and $\Omega$ are real analytic and $\Omega(0,0)\neq0$.
By using \eqref{eq-Rxy} and the fact that $\Omega(0,0)\neq0$, we may let
$$R(0,y)=2y^2+\sum_{n=3}^{\infty}c_n^{*}y^n~\mbox{ and }~\Omega(0,y)=\Omega(0,0)+\sum_{n=1}^{\infty}d_n^{*}y^n
$$
so that \eqref{eq-Rxy-extra} takes the form
$$
2y^2+\sum_{n=3}^{\infty}c_n^{*}y^n=\big(y^2+A_1(0)y+A_0(0)\big)\left (\Omega(0,0)+\sum_{n=1}^{\infty}d_n^{*}y^n\right ).
$$
Comparing the coefficients of $y^k$ ($k=0,1,2$) gives
$$A_0(0)\Omega(0,0)=0,\,A_1(0)\Omega(0,0)+A_0(0)d_1^{*}=0
$$
and
$$A_0(0)d_2^{*}+A_1(0)d_1^{*}+\Omega(0,0)=2
$$
from which it follows that $A_0(0)=0=A_1(0)$ and thus, $\Omega(0,0)=2$.

Clearly, there exists  $\delta_{3}'''\in (0, \delta_3'']$ such that for $z$ with $|z|<\delta_{3}'''$,
we have
\be\label{eq2.6}
1<\Omega(x,y)<4
\ee
and
\be\label{eq2.5}
\left |\frac{b_{2m+1} +\sum_{k=2m+2}^{\infty}b_{k}z^{k-2m-1}}{1+\sum_{k=2}^{\infty}a_kz^k}\right |>\frac{|b_{2m+1}|}{2}.
\ee
Let us now introduce
$$ \mathcal {G}(x,y)=y^2+A_1(x)y+A_0(x).
$$
Then the inequalities \eqref{eq-Rxy-extra}, \eqref{eq2.6}, \eqref{eq-Rxy},
\eqref{eq2.4} and \eqref{eq2.5} show that
\be\label{eq2.5-e1} \mathcal {G}(x,y)\geq
\frac{|b_{2m+1}|}{8}|z|^{2m+1} ~\mbox{ for $z$ with
$0<|z|<\delta_{3}'''$}.
\ee
We have seen that $A_0(0)=0$ and $A_1(0)=0$.
  Hence the function $A_1$ has the following power series expansion:
$$A_1(x)=c_{k_1}x^{k_1}+\sum_{k=k_1+1}^{\infty}c_kx^k,
$$
where $k_1\geq 1$ and $c_{k_1}\neq 0$. Again, since
$$\mathcal{G}(x,y)=\left (y+\frac{A_1(x)}{2}\right )^2+\left (A_0(x)-\frac{A_1^2(x)}{4}\right ),
$$
it follows  from the last relation and \eqref{eq2.5-e1} that for $0<\sqrt{x^2+\frac{A_1^2(x)}{4}}<\delta_{3}'''$,
$$ \mathcal {G}\Big(x, -\frac{A_1(x)}{2}\Big)=A_0(x)-\frac{A_1^2(x)}{4}\geq
\frac{|b_{2m+1}|}{8}|x|^{2m+1}.
$$
Consequently, we have the following series expression:
$$ A_0(x)-\frac{A_1^2(x)}{4}=d_{2m_0}x^{2m_0}+\sum_{k=2m_0+1}^{\infty}d_kx_k,
$$
where $m_0\geq 1$, $d_{2m_0}\neq 0$ and $2m_0< 2m+1$. 
Note that the series must begin with an even power since the function is positive for small $x$.

Using  similar arguments as  in the proof of \cite[Lemma 1]{CW}, we
find that there exist $\delta_{3}'''$ with $0<\delta_{3}''''<\delta_{3}'''$ and $C_0>0$
such that
$$\mathcal {G}(x,y)\geq C_0|z|^{2m_0}, ~\mbox{ for $z$ with $0<|z|<\delta_{3}''''$}
$$
and then, by \eqref{eq2.6},
\be\label{eq2.19}
R(x,y)\geq C_0|z|^{2m_0}.\ee Let
$$M_0=\sup_{|z|\leq \delta_3'''}\left \{\frac{1}{|1+\sum_{k=2}^{\infty}a_kz^k|}\right \}
$$
and
$$M_1=\sup_{|z|\leq \delta_3'''}\left \{\left |\frac{b_{2m+1}+\sum_{k=2m+2}^{\infty}b_{k}
z^{k-2m-1}}{1+\sum_{k=2}^{\infty}a_kz^k}\right |\right \}.
$$
Since $2m_0<2m+1$, we can choose $0<\delta_3\leq \delta_{3}''''$
such that
$$ C_0|z|^{2m_0}\geq (M_0+M_1)|z|^{2m+1} ~\mbox{ for $z$ with $0<|z|<\delta_3$}.
$$
Then \eqref{eq2.19} shows
\beqq
R(x,y)
&\geq& (M_0+M_1)|z|^{2m+1}\\
&\geq& \frac{|z|^{2m+1}}{|1+\sum_{k=2}^{\infty}a_kz^k|}+
\left|\frac{b_{2m+1}z^{2m+1}+\sum_{k=2m+2}^{\infty}b_{k}
z^{k}}{1+\sum_{k=2}^{\infty}a_kz^k}\right| =\frac{|z|^{2m+1}+|g'(z)|}{|h'(z)|},
 \eeqq
which by the definition of $R(x,y)$ given by \eqref{eq-Rxy} implies that
\beq\label{qw-4}
\big(|h'(z)|+|g'(z)|+|z|^{n_3}\big)(1-|z|^2)<1
\eeq
for $n_3=2m+1$. Thus, \eqref{eq-ass1a} holds with $n=2m+1$ and $\delta =\delta_3$.

\bsca\label{qw-8} $|b_2|=1$.\esca

Clearly, $a_2=0$ in this case.
By considering the function $h+e^{i\theta}\overline{g}$, if needed, we assume that $b_2=1$, and thus,
$h'(z)$ and $g'(z)$ take the form
$$h'(z)=1+a_3z^3+\sum_{k=4}^{\infty}a_kz^k\;\;\mbox{and}\;\;
g'(z)=z^2+b_3z^3+\sum_{k=4}^{\infty}b_kz^k.
$$
Further, we assume that $a_3\geq0$.

Obviously, $h'(z)$ is non-vanishing and analytic in
$|z|<\delta_0'$  for some $0<\delta_0'<\delta_0$ and therefore, there exists an analytic function $\Psi_0$ such that $h'(z)=\Psi_0^2(z)$ in
$|z|<\delta_0'$. As $h'(0)=1$, $h''(0)=0$ and $h'''(0)=0$, we find that
$\Psi_0 (0)=1$,  $\Psi_0'(0)=0$, $\Psi_0'' (0)=0$ and $\Psi_0'''(0)=h^{(4)}(0)/2\, =3!a_3/2$. Thus, we have the following series expansion for
$\Psi_0$:
$$ \Psi_0 (z)=1 +\frac{a_3}{2}z^3 +\cdots ~\mbox{ for $|z|<\delta_0'$}.
$$
Clearly, the last relation yields
$$|h'(z)|=\Psi_0 (z)\overline{\Psi_0(z)}=\left (1+\frac{a_3}{2}z^3 +\cdots\right )\overline{\left ( 1 +\frac{a_3}{2}z^3 +\cdots\right)}=1+{\rm Re\,}(a_3 z^3)+ O_4(x,y).
$$

A similar procedure for the function $g(z)/z^2$ gives
$$\left |\frac{g'(z)}{z^2}\right |=1+{\rm Re\,}(b_3z)+O_2(x,y)
$$
for $0<|z|<\delta_0''$ with $0<\delta_0''\leq \delta_0'$, which shows that
$$|g'(z)|=|z|^2(1+{\rm Re\,}(b_3z)+O_2(x,y))  ~\mbox{ for $0<|z|<\delta_0''$.}
$$

It follows from the assumption \eqref{eq-ass1} that
\beq\label{eq2.12}
|h'(z)|+|g'(z)|&=&1+|z|^2+{\rm Re\,}(a_3
z^3+b_3|z|^2z)+
O_4(x,y)\\
&<&\nonumber\frac{1}{1-|z|^2}
=1+|z|^2+|z|^4+\sum_{k=3}^{\infty}|z|^{2k}
\eeq 
for $0<|z|<\delta_0''$.

\bcl \label{cl-2} $a_3=b_3=0.$\ecl

Suppose on the contrary that either $a_3\neq 0$ or $b_3\neq 0$. Letting
$z=r$ in \eqref{eq2.12} leads to
$$a_3+{\rm Re\,}(b_3)=0.
$$

If $a_3=|b_3|$, that is, $a_3=-b_3\neq 0$, then for $z=re^{i\theta}$
with $0<r<\delta_0$ and $\cos3\theta-\cos \theta \neq 0$, \eqref{eq2.12} yields
$$a_3(\cos3\theta-\cos \theta)r^3+O_4(x,y)< r^4+\sum_{k=3}^{\infty}r^{2k}.
$$
This obvious contradiction shows that $a_3=0$ and thus, $b_3=0$.

If $a_3\neq|b_3|$, assume first
$|b_3|>a_3$. Thus, for $z=re^{i\theta_0}$ with $0<r<\delta_0$ and
$\theta_0=-\arg b_3$, we have
$${\rm Re\,}(a_3 z^3+b_3|z|^2z)=|b_3|r^3+{\rm Re\,}(a_3 z^3)\geq (|b_3|-a_3)r^3,
$$
and then, we infer from \eqref{eq2.12} that
$$(|b_3|-a_3)r^3+O_4(x,y)< r^4+\sum_{k=3}^{\infty}r^{2k},
$$
which is again a contradiction. If $a_3<|b_3|$, we obtain a similar contradiction for $z=r$.
The proof of Claim \ref{cl-2} is finished.
\medskip

Now, by Claim \ref{cl-2}, it is easy to show that
$$|h'(z)|=\Big|1+a_4z^4+\sum_{k=5}^{\infty}a_kz^k\Big|=1+{\rm Re\,}(a_4z^4)+O_5(x,y)
$$
and
$$|g'(z)|=|z|^2\big(1+{\rm Re\,}(b_4z^2)+O_3(x,y)\big).
$$

Let
$$ \mathcal {H}(x,y)=\big(1-(x^2+y^2)\big)^{-1}-\big(|h'(x+iy)|+|g'(x+iy)|\big).
$$
By the assumption  \eqref{eq-ass1}, we observe that
$$\mathcal {H}(x,y)>0\;\;\mbox{ for}\;\; 0<|z|<\delta_0
$$
and, by the representation of $|h'(z)|$, $|g'(z)|$ and \eqref{eq-ass1}, we have the inequality
\be\nonumber
{\rm Re\,}(a_4z^4)+|z|^2{\rm Re\,}(b_4z^2)\leq |z|^4
\ee 
for $0<|z|<\delta_0'''$ with $0<\delta_0'''\leq \delta_0$.
By using a rotation and without loss of generalization, we assume $a_4\geq 0$. Consequently, the last
relation is equivalent to
\be\label{eq2.24}
a_4{\rm Re\,}(z^4)+|z|^2{\rm Re\,}(b_4z^2)\leq |z|^4\;\;\mbox{for}\;\; 0<|z|<\delta_0'''.
\ee
Let $z=x+iy$ with $0<x^2+y^2<\delta_0'''^2$ and $b_4=a+ib$. Then \eqref{eq2.24} is equivalent to
$$(a_4+a)x^4+(a_4-a)y^4-2bxy(x^2+y^2)-6a_4x^2y^2\leq x^4+y^4+2x^2y^2,
$$
which implies that $a_4+a\leq 1$ and $a_4-a\leq 1$, i.e. $a_4+|a|\leq 1$. Thus, we may rewrite the
last relation as
\be\label{eq2.27}
(1-a_4-a)x^4+(1-a_4+a)y^4+2bxy(x^2+y^2)+(6a_4+2)x^2y^2\geq 0
\ee
for $0<x^2+y^2<\delta_0'''^2$. Now, without loss of generalization, we assume $a\geq0$.

If $a_4+a=1$, then \eqref{eq2.27} becomes
\be \label{eq2.25}
(1-a_4+a)y^4+2bxy(x^2+y^2)+(6a_4+2)x^2y^2\geq 0
\ee
for $0<x^2+y^2<\delta_0'''^2$.  Now, we prove $b=0$. Suppose not. Then $b\neq0$  and
\be \label{eq2.31}
(1-a_4+a)y^4+(6a_4+2)x^2y^2\geq |2bxy|(x^2+y^2)
\ee
with $0<x^2+y^2<\delta_0'''^2$.
If $|y|<\min\{\frac{2|b|}{1-a_4+a},\frac{|b|}{3a_4+1}\}|x|$ with $0<x^2+y^2<\delta_0'''^2$, then  we have
$$
(1-a_4+a)y^4<2|b||xy^3|\;\;\mbox{and}\;\;(6a_4+2)x^2y^2<2|b||x^3y|
$$
which contradicts with the inequality \eqref{eq2.31}.  Hence, we must have $b=0$, and thus
$$\mathcal {H}(x,y)=(1-a_4+a)y^4+{\color{red}{(6a_4+2)}}x^2y^2+O_5(x,y)>0
$$
for $0<x^2+y^2<\delta_0^2$. It is easy to verify that there are a $\delta_4'$ with
$0<\delta_4'<\delta_0$ and an integer $n_4'\geq 2$  such that
\be\label{eq2.28}
\mathcal {H}(x,y)>(x^2+y^2)^{n_4'},
\ee
where $0<x^2+y^2<(\delta_4')^2$ and $\mathcal {H}(x,0)=\sum_{n=2(n_4'-1)}^{\infty}\alpha_n x^n$.

Now we discuss the case $a_4+a<1$ such that \eqref{eq2.27} holds
for $0<x^2+y^2<\delta_0'''^2$. We will now prove that if $b\neq 0$, then
$(1-a_4-a)(3a_4+1)\geq b^2$. If not there exist $x$ and
$y$ with $0<x^2+y^2<\delta_0'''^2$ such that
$$\frac{1-a_4-a}{b}|x|<|y|<\frac{b}{3a_4+1}|x|
$$
from which we obtain that
$$\left \{ \begin{array}{l}
(1-a_4-a)x^4<|bxy^3|,\\
(1-a_4+a)y^4<|bxy^3|,\\
(6a_4+2)x^2y^2<|2bx^3y|,
\end{array} \right.
$$
which contradicts the inequality \eqref{eq2.27}. Hence
$$(1-a_4-a)(3a_4+1)\geq b^2.
$$

For the case $b\neq 0$ and $a\neq0$,  by using the inequality $(1-a_4-a)(3a_4+1)\geq b^2$, we have
\beqq
&&(1-a_4-a)x^4+(1-a_4+a)y^4+(6a_4+2)x^2y^2\\
 &>&
(1-a_4-a)x^4+(3a_4+1+\varepsilon_0)x^2y^2+(1-a_4+a)y^4+(3a_4+1-\varepsilon_0)x^2y^2\\
&\geq& 2\sqrt{1-a_4-a}\sqrt{3a_4+1+\varepsilon_0}|x y|x^2+2\sqrt{1-a_4+a}\sqrt{3a_4+1-\varepsilon_0}|x y|y^2\\
&=&2\sqrt{1-a_4-a}\sqrt{3a_4+1+\varepsilon_0}| x y|(x^2+y^2)\\
&>&2|b x y|(x^2+y^2)
\eeqq
with $\varepsilon_0=\frac{a(3a_4+1)}{1-a_4}$,
which implies that there exists  $\delta_4''$ with $0<\delta_4''<\delta_0'''$, such that
\beq\label{eq2.30}
&&\mathcal {H}(x,y)\\&=&\nonumber(1-a_4-a)x^4+(1-a_4+a)y^4+2bxy(x^2+y^2)+(6a_4+2)x^2y^2+O_5(x,y)\\
&\geq&\nonumber (1-\frac{b}{B})\Big((1-a_4-a)x^4+(1-a_4+a)y^4+(6a_4+2)x^2y^2\Big)+O_5(x,y)\\
&>& \nonumber(x^2+y^2)^3
\eeq
for $0<x^2+y^2<(\delta_4'')^2$, $B=\sqrt{1-a_4-a}\sqrt{3a_4+1+\varepsilon_0}$.

We claim that if $b\neq 0$ and $a=0$, then $(1-a_4)(3a_4+1)>b^2$. Otherwise $(1-a_4)(3a_4+1)=b^2$.
Let $\varepsilon_1=\sqrt{(3a_4+1)^2-(1-a_4)^2}$. Then
\beq\label{2.30}
&&(1-a_4)x^4+(1-a_4)y^4+2bxy(x^2+y^2)+(6a_4+2)x^2y^2\\
\nonumber &\geq&  \sqrt{1-a_4}\Big(\sqrt{3a_4+1+\varepsilon_1}+ \sqrt{3a_4+1-\varepsilon_1}\Big)|xy|(x^2+y^2),
\eeq
and the equality holds when $|y|=\frac{1-a_4}{3a_4+1+\varepsilon_1}|x|$.
But this is a contradiction since $\sqrt{3a_4+1+\varepsilon_1}+ \sqrt{3a_4+1-\varepsilon_1}<2\sqrt{3a_4+1}$.
For the case $b\neq 0$ and $a=0$, by using the inequality $(1-a_4)(3a_4+1)>b^2$ and the similar arguments as that of
the case $b\neq 0$ and $a\neq 0$, we obtain that
there exists a $\delta_4'''$
with $0<\delta_4'''<\delta_0'''$ such that
\be\label{2.29}
\mathcal {H}(x,y)>(x^2+y^2)^3
\ee
for $0<x^2+y^2<(\delta_4''')^2$.

If $b=0$, then
$$\mathcal {H}(x,y)=(1-a_4-a)x^4+(1-a_4+a)y^4+(6a_4+2)x^2y^2+O_5(x,y)>0
$$
for $0<x^2+y^2<\delta_0^2$. It is easy to verify that there is a $\delta_4''''$
with $0<\delta_4''''<\delta_0$ such that
\be\label{eq2.29}
\mathcal {H}(x,y)>(x^2+y^2)^3
\ee
for $0<x^2+y^2<(\delta_4''')^2$.

Therefore, for $n_4=\min\{n_4', 6\}$ and $\delta_4=\min \{\delta_4',
\delta_4'', \delta_4''', \delta_4''''\}$, we deduce that
\be\label{qw-5}
(h'(z)+g'(z)+|z|^{n_4})(1-|z|^2)<1 ~\mbox{ for $0<|z|<\delta_4$.}
\ee


\medskip

Finally, we let
$\delta=\min\{\delta_1, \delta_2,\delta_3', \delta_3, \delta_4\}$, $n= \max\{n_1, n_2, n_3, n_4\}$,
and observe that the inequalities \eqref{qw-1},  \eqref{qw-2}, \eqref{qw-3}, \eqref{qw-4}
and \eqref{qw-5} show that
$$\big(|h'(z)|+|g'(z)|+|z|^n\big)(1-|z|^2)<1 ~\mbox{for all $z$ with $0<|z|<\delta$.}
$$
Thus, \eqref{eq-ass1a} holds  and the proof of Lemma \ref{lem3.2} is complete.
\epf

\blem\label{lem3.3} Suppose
$f=h+\overline{g}\in\widetilde{\mathscr{B}}_{\mathcal {H},0,1}$ $($resp.
$f=h+\overline{g}\in\widetilde{\mathscr{B}}_{\mathcal {H},1}$ $)$
such
 that
 \bee
\item there is a $z_0\in \mathbb{D}$ with
$|h'(z_0)|(1-|z_0|^2)=1$ or $|g'(z_0)|(1-|z_0|^2)=1$; and
\item there exists a $\delta>0$ such that
$$\big(|h'(z)|+|g'(z)|\big)(1-|z|^2)<1  ~\mbox{ for $0<|z-z_0|<\delta$.}
$$
\eee
Then there are a positive integer $n$ and a $\delta'\in (0,
\delta]$ such that
$$\Big(|h'(z)|+|g'(z)|+\Big|\frac{z-z_0}{1-\overline{z_0}z}\Big|^n\Big)\big(1-|z|^2\big)<1
~\mbox{ for $0<|z|<\delta'$.}
$$
\elem

\bpf The proof follows easily from Lemma \ref{lem3.2} and the similar
argument as in the proof of \cite[Lemma 3]{CW}.
\epf

\setcounter{case}{0} \setcounter{subcase}{0}

\blem\label{lem3.4} Suppose that $f=h+\overline{g}\in
\widetilde{\mathscr{B}}_{\mathcal {H},1}$ satisfies the following conditions:
\bee
\item $h'(0)\neq 0$, $g'(0)\neq 0$ and $|h'(0)|+|g'(0)|=1$;
\item there is a $\delta_0>0$ such
that
\be\label{eq-ass3}
\big(|h'(z)|+|g'(z)|\big)(1-|z|^2)<1~\mbox{ for  $0<|z|<\delta_0$.}
\ee
\eee
Then there are a positive integer $n$ and
$\delta \in (0,  \delta_0]$ such that
$$\big(|h'(z)|+|g'(z)|+|z|^n\big)(1-|z|^2)<1~\mbox{ for $0<|z|<\delta$.}
$$
\elem

\bpf  By considering the function
$e^{i\theta_1}h+e^{i\theta_2}\overline{g}$, if needed,
 we may assume that $h'(0)\in (0,1)$ and
$g'(0)\in (0,1)$. Then there is a $\delta'\in (0, \delta_0]$ such
that $h'(z)\neq 0$ and $g'(z)\neq 0$ in $\mathbb{D}_{\delta'}$. We
assume that
$$h'(z)=a_0+a_1z+\sum_{k=2}^{\infty}a_kz^k \;\,\mbox{ and}\;\,
g'(z)=b_0+b_1z+\sum_{k=2}^{\infty}b_kz^k .
$$

Since $h$, $g$ and $h+g$ are (analytic) Bloch functions, it follows
from the similar reasoning as in the proof of \cite[Lemma 2]{CW}
that $a_1=b_1=0$, $|a_2+b_2|\leq 1$. Without loss of generality, we
assume that $0\leq a_2+b_2\leq1$.  Let
$$ \mathcal {H}(x,y)=\big(1-(x^2+y^2)\big)^{-1}-\big(|h'(x+iy)|+|g'(x+iy)|\big).
$$
By the assumption  \eqref{eq-ass3}, we see that
\be\mathcal {H}(x,y)>0\;\;\mbox{ for}\;\; 0<|z|<\delta_0.
\ee
It is easy to verify that $\mathcal {H}(0,y)=\sum_{n=2}^{\infty}c_ny^n$ with $c_2>0$.  It follows from
Theorem \Ref{Lem2.0} that
$$\mathcal {H}(x,y)=\big(y^2+A_3(x)y+A_4(x)\big)F (x,y),
$$
where $A_3$, $A_4$ and $F$ are real analytic functions, and $F(0,0)\neq 0$
(Actually $F(0,0)=c_2$). By using \cite[Lemma 1]{CW}, there are an $n_0$ and
a $0<\delta<\delta_0$ so that for $0<x^2+y^2<\delta^2$,
$$y^2+A_3(x) y+A_4(x)>(x^2+y^2)^{n_0}.
$$
Since $F(0,0)=c_2>0$, a possible smaller choice of $\delta$ yields that if $0<x^2+y^2<\delta^2$, then
$$\mathcal {H}(x,y)>(x^2+y^2)^{n_0},
$$
which implies that
$$\big(|h'(z)|+|g'(z)|+|z|^{n}\big)\big(1-|z|^2\big)<1
$$
for $0<|z|<\delta$ and $n=2n_0$.
\epf

\blem\label{lem3.5} Suppose $f=h+\overline{g}\in
\widetilde{\mathscr{B}}_{\mathcal {H},0,1}$ $($resp.
$f=h+\overline{g}\in \widetilde{\mathscr{B}}_{\mathcal {H},1}$ $)$
satisfies the following:
\bee
\item There is a $z_0\in \mathbb{D}$ such that $h'(z_0)\neq 0$, $g'(z_0)\neq 0$
and $\big(|h'(z_0)|+|g'(z_0)|\big)(1-|z_0|^2)=1$;
\item There exists a $\delta>0$ so that for all $z$ with $0<|z-z_0|<\delta$,
$$\big(|h'(z)|+|g'(z)|\big)\big(1-|z|^2\big)<1.
$$
\eee
Then there are a positive integer $n$ and
$\delta'$ with $0<\delta'\leq \delta$ such that
$$\left (|h'(z)|+|g'(z)|+\left |\frac{z-z_0}{1-\overline{z_0}z}\right |^n\right )\left(1-|z|^2\right )<1 ~\mbox{ for  $0<|z|<\delta'$.}
$$
\elem \bpf The proof of Lemma \ref{lem3.5}  follows easily from Lemma \ref{lem3.4} and the
similar reasoning as in the proof of \cite[Lemma 3]{CW}.
\epf

\subsection{The proof of  Theorem \ref{thm2.0}}

 The proof easily follows from Lemmas \ref{lem3.3},
\ref{lem3.5} and the similar reasoning as in the proof of
\cite[Theorem 2]{CW}.

\section{Support points }\label{sec-3}
Let $\mathcal {L}$ be a continuous linear functional of $\mathcal
{A}(\mathbb{D})$ into $\mathbb{C}$. By  \cite[Theorem 9.3]{Du2}, we
know that there must be a one-to-one correspondence between
$\mathcal {L}$ and the set of sequences $\{A_k\}$  of complex numbers with
$$\limsup_{k\rightarrow\infty}|A_k|^{\frac{1}{k}}<1,
$$
and for analytic functions $h$, if $h(z)=\sum_{k=0}^{\infty}a_kz_k$, then
$\mathcal {L}(h)=\sum_{k=0}^{\infty}A_ka_k$. For the setting of
harmonic mappings, we have the following analog of it.

\blem\label{lem2.1} $($ \cite[pp.~131]{SS}$)$ Suppose that $\mathcal {L}$ is
a continuous linear functional of $\mathcal {H}(\mathbb{D})$. Then
there are two sequences $\{A_k\}$ and $\{B_k\}$ such that
\bee
\item[ {\rm (1)}] $\limsup_{k\rightarrow\infty}|A_k|^{\frac{1}{k}}<1$,
$\limsup_{k\rightarrow\infty}|B_k|^{\frac{1}{k}}<1$, and

\item[ {\rm (2)}] for $f\in \mathcal {H}(\mathbb{D})$, if
$$f(z)=h(z)+\overline{g}(z)=\sum_{k=0}^{\infty}a_kz^k+\sum_{k=0}^{\infty}\overline{b}_k\bar{z}^k,
$$
then
$$\mathcal {L}(f)=\sum_{k=0}^{\infty}A_k a_k+\sum_{k=0}^{\infty}\overline{B}_k \overline{b}_k.
$$
\eee
\elem

From this lemma, we observe that there is a one-to-one correspondence
between continuous linear functionals $\mathcal {L}$ and pair of sequences
$\{A_k\}$ and $\{B_k\}$ of complex numbers with
$$\limsup_{k\rightarrow\infty}|A_k|^{\frac{1}{k}}<1\;\;\mbox{and}\;\;\limsup_{k\rightarrow\infty}|B_k|^{\frac{1}{k}}<1,
$$
respectively.

Now, we introduce some lemmas which are useful in the proof of our
main result of this section.

\blem\label{lem2.2}
Suppose that $\mathcal {L}$ is a continuous linear functional and that $f\in \mathcal
{H}(\mathbb{D})$. For $\varepsilon \in (0,1]$, we define
$f_\varepsilon\in\mathcal {H}(\mathbb{D})$ by
$$f_\varepsilon(z)=f((1-\varepsilon)z)
$$
in $\mathbb{D}$. Then there exists a constant $K>0$ such that $|\mathcal {L}(f_\varepsilon-f)|\leq \varepsilon K.$
\elem

\bpf The proof follows  from Lemma \ref{lem2.1} and the similar argument as in the proof of
\cite[Lemma 3]{Bo}.
\epf

\begin{Lem}\label{lem2.3}$($\cite[Lemma 4]{Bo}$)$ Suppose $M\geq 0$. Then
there exist numbers $\varepsilon_1$ and $R\in (0,1)$ such that for
all $\varepsilon$ with $0<\varepsilon\leq\varepsilon_1$ and all $z$ with
$R\leq |z|<1$,
$$\frac{1}{1-(1-\varepsilon)^2|z|^2}+\frac{\varepsilon M}{1-|z|^2}\leq
\frac{1}{1-|z|^2}.
$$
\end{Lem}

In order to state and prove the next lemma, we introduce the following notations.

For $f=h+\overline{g}\in \mathcal {H}(\mathbb{D})$, denote
$$M(f)=\sup_{z\in
\mathbb{D}}\Big\{\big(|h(z)|+|g(z)|\big)(1-|z|^2)\Big\},
$$
$$\Gamma(f)=\{z\in \mathbb{D}: \big(|h(z)|+|g(z)|\big)(1-|z|^2)=1\},
$$
and $K_{\mathcal {H},1}=\{f: M(f)\leq 1\}.$

 \blem\label{lem2.4} Suppose that
$\mathcal {L}$ $(\mathcal {L}\not\equiv 0)$ is a continuous linear
functional of $\mathcal {H}(\mathbb{D})$ and that $\sup_{f^*\in
K_{\mathcal {H},1}}\Big\{{\rm Re\,}
\{\mathcal {L}(f^*)\}\Big\}= {\rm Re\,}
\{\mathcal {L}(f)\}$. Then $\Gamma (f)\neq\emptyset$. \elem

\bpf Suppose on the contrary that $\Gamma(f)=\emptyset$. This means that if
$f=h+\overline{g}$, then for $z\in\mathbb{D}$,
$$\big(|h(z)|+|g(z)|\big)(1-|z|^2)<1.
$$

Assume that the functional $\mathcal {L}$ corresponds to the sequences $\{A_k\}$ and
$\{B_k\}$. Since $\mathcal {L}\not\equiv 0$, we may assume that $A_{k_0}\neq 0$.

Let $H(z)=\frac{2K}{A_{k_0}}z^{k_0}$, where $K$ is the same as in
Lemma \ref{lem2.2} with respect to $\mathcal {L}$ and $f$. Clearly, it is
easy to verify that
$$M=M(H)<\infty,\;\;{\rm Re\,}{\mathcal L}\,(H)=2K\;\,\mbox{and}\;\, H\not\equiv0.
$$
Define $\widetilde{f}$ by
$$ \widetilde{f}(z)= f\big((1-\varepsilon)z)+\varepsilon H((1-\varepsilon)z\big),
$$
where $\varepsilon\in (0,\varepsilon_1)$, and $\varepsilon_1$ is
the same as in Lemma \Ref{lem2.3} with the constant  $M(H)$ in place of $M$.
Obviously, there exists a $\delta>0$ such that
\be\label{eq3.2}
|h(z)|+|g(z)|+\delta\leq \frac{1}{1-|z|^2}
\ee
in $\mathbb{D}(R)$, where  $R$ is the same as in Lemma \Ref{lem2.3} with the constant
$M(H)$ in place of $M$.

Finally, we choose $\varepsilon_2\in (0, \varepsilon_1)$ such that
\be\label{eq3.3}
{\rm Re\,}{\mathcal L}(H((1-\varepsilon)z))\geq \frac{3}{2}K
\ee
and
\be\label{eq3.4}
\varepsilon|H((1-\varepsilon)z)|\leq \delta
\ee
in $\mathbb{D}(R)$ for $\varepsilon\in(0,\varepsilon_2].$
Then  \eqref{eq3.2} and \eqref{eq3.4} imply that for $z\in \mathbb{D}(R)$ and
$\varepsilon\in(0,\varepsilon_2]$,
\beq\label{eq3.5}
|h((1-\varepsilon)z)|+|g((1-\varepsilon)z)|+\varepsilon|H((1-\varepsilon)z)|
&\leq&\nonumber|h((1-\varepsilon)z)|+|g((1-\varepsilon)z)|+\delta\\
&\leq&\nonumber\frac{1}{1-(1-\varepsilon)^2|z|^2}\\
&<&\frac{1}{1-|z|^2}.
\eeq
We infer from Lemma \Ref{lem2.3} that for $z$ with $R\leq|z|<1$ and
$\varepsilon\in(0,\varepsilon_2]$,
\beq\label{eq3.6}
|h((1-\varepsilon)z)|+|g((1-\varepsilon)z)|+\varepsilon|H((1-\varepsilon)z)|
&\leq&\nonumber\frac{1}{1-(1-\varepsilon)^2|z|^2}+\frac{\varepsilon
M}{1-|z|^2}\\
&\leq&\frac{1}{1-|z|^2}.
\eeq
Therefore we know from \eqref{eq3.5} and \eqref{eq3.6} that $\widetilde{f}\in K_{\mathcal
{H},1}$ and thus,
\beqq
{\rm Re\,}\mathcal {L}\big(\widetilde{f})-{\rm Re\,}\mathcal {L}(f) &=&{\rm Re\,}\mathcal
{L}\big(f((1-\varepsilon)z)-f(z)\big)+\varepsilon {\rm Re\,}\mathcal
{L}\big(H((1-\varepsilon)z)\big)\\
&\geq& \varepsilon {\rm Re\,}\mathcal
{L}\big(H((1-\varepsilon)z)\big)-\varepsilon K ~\mbox{ (by Lemma \ref{lem2.2}) }\\
&\geq& \frac{3}{2}\varepsilon K - \varepsilon K=\frac{K}{2}\varepsilon\\
&>&0,
\eeqq
which is a contradiction, because
$$\ds\sup_{f^*\in K_{\mathcal {H},1}}\big\{{\rm Re\,}\{\mathcal {L}(f^*)\}\big\}= {\rm Re\,} \{\mathcal {L}(f)\}.
$$
The proof of the lemma is complete.
 \epf

Now, we are ready to state and prove our main result of this section,
which is a characterization of support points in the unit ball of harmonic Bloch spaces in
$\mathbb{D}$.

\bthm\label{thm2.1} We have
\bee
\item[ {\rm (1)}] A function $f_0\in \mathscr{B}_{\mathcal {H}, 1}$
is a support point of $\mathscr{B}_{\mathcal {H}, 1}$ if and only if $f_0$ is
a convex combination of a unimodular constant $u$ and a support
point $f$ of $\widetilde{\mathscr{B}}_{\mathcal {H}, 1}$, i.e., there
are constants $\lambda_1$, $\lambda_2 \in [0,1]$ with
$\lambda_1+\lambda_2=1$ such that
$$f_0=\lambda_1u+\lambda_2f.
$$



\item[ {\rm (2)}] A function $f\in \widetilde{\mathscr{B}}_{\mathcal {H},
1}$ is a support point of $\widetilde{\mathscr{B}}_{\mathcal {H}, 1}$ if
and only if $\Lambda_f\neq \emptyset$.
\eee\ethm

\bpf $(1)$: The proof of the first part easily follows from the similar reasoning as in the proof  of \cite[Corollary
$2$]{CW}.

 $(2)$:  For the sufficiency of the second part of the theorem, we assume that
$\Lambda_f\neq \emptyset$ for
$f=h+\overline{g}\in \widetilde{\mathscr{B}}_{\mathcal {H}, 1}$, and let $z_0\in \Lambda_f$. Then
\be\label{eq3.7}
|h'(z_0)|+|g'(z_0)|=\frac{1}{1-|z_0|^2}=\sup_{f_1=h_1+\overline{g_1}\in
\widetilde{\mathscr{B}}_{\mathcal {H},
1}}\big(|h_1'(z_0)|+|g_1'(z_0)|\big).
\ee
Let $\theta_0$ be such that
\be\label{eq3.8}
|h'(z_0)|+|g'(z_0)|=|h'(z_0)+e^{i\theta_0}\overline{g'(z_0)}|.
\ee
Then \eqref{eq3.7} and \eqref{eq3.8} yield that
\be\label{eq3.9}
|h'(z_0)+e^{i\theta_0}\overline{g'(z_0)}|
=\sup_{f_1=h_1+\overline{g_1}\in \widetilde{\mathscr{B}}_{\mathcal {H},
1}}\big(|h_1'(z_0)+e^{i\theta_0}\overline{g_1'(z_0)}|\big).
\ee
For $f_1=h_1+\overline{g_1}\in \widetilde{\mathscr{B}}_{\mathcal {H}, 1}$,
we define
$$\mathcal {L}(f_1)=\big(\overline{h'(z_0)}+e^{-i\theta_0}\overline{g'(z_0)}\big)\big(h_1'(z_0)+e^{i\theta_0}\overline{g_1'(z_0)}\big).
$$
Then $\mathcal {L}$ is a continuous linear functional, and by \eqref{eq3.9}, we have
$${\rm Re\,}\mathcal {L}(f_1)\leq |h'(z_0)+e^{i\theta_0}\overline{g'(z_0)}|^2={\rm Re\,} \mathcal {L}(f).
$$
It follows that $f$ is a support point of $\widetilde{\mathscr{B}}_{\mathcal {H},1}$.

Conversely, assume that $f$ is a support point of
$\widetilde{\mathscr{B}}_{\mathcal {H},1}$. Then there exists a
continuous linear functional $\widetilde{\mathcal {L}}$ such that
$\widetilde{\mathcal {L}}$ is not constant on
$\widetilde{\mathscr{B}}_{\mathcal {H},1}$ and ${\rm Re\,}
\widetilde{\mathcal {L}}(f_1)\leq {\rm Re\,} \widetilde{\mathcal {L}}(f)$ for
each $f_1=h_1+\overline{g_1}\in \widetilde{\mathscr{B}}_{\mathcal
{H},1}$. We assume that
$$\widetilde{\mathcal {L}}(f_1)=\sum_{k=1}^{\infty}A_k a_k+\sum_{k=1}^{\infty}\overline{B_k}\,\overline{b_k}
$$
for
$$f_1(z)=h_1(z)+\overline{g_1 (z)}=\sum_{k=1}^{\infty}a_kz^k+\sum_{k=1}^{\infty}\overline{b_k}\bar{z}^k.
$$
Define $C_k=\frac{A_{k+1}}{k+1}$ and $D_k=\frac{B_{k+1}}{k+1}$. Obviously,
$$\limsup_{k\rightarrow\infty}|C_k|^{\frac{1}{k}}<1\;\;\mbox{and}\;\;\limsup_{k\rightarrow\infty}|D_k|^{\frac{1}{k}}<1.
$$
Consider the continuous linear functional
$$L(f_1)=\sum_{k=1}^{\infty}C_k
a_k+\sum_{k=1}^{\infty}\overline{D}_k\overline{ b}_k.
$$
Then
$\widetilde{\mathcal {L}}(h_1+\overline{g_1})=L(h_1'+\overline{g_1'})$ and ${\rm Re\,} L(h_1'+\overline{g_1'})\leq {\rm Re\,} L(h'+\overline{g'}).
$
By Lemma \ref{lem2.4}, we see that
$\Lambda_f=\Gamma(h'+\overline{g'})\neq\emptyset$, which implies that
Theorem \ref{thm2.1} is true.\epf

\subsection*{\bf Acknowledgements}
The  work of the second author is
supported by Mathematical Research Impact Centric Support of DST, India  (MTR/2017/000367).
The third author is supported by National Natural
Science Foundation of China (No. 11501159).

\subsection*{Conflict of Interests}
The authors declare that there is no conflict of interests regarding the publication of this paper.


\begin{thebibliography}{99}

\bibitem{Am} { Y. Abu-Muhanna}, On extreme points of subordination
families. \textit{Proc. Amer. Math. Soc.} \textbf{87} (1983), 439--443.

\bibitem{Am1} { Y. Abu-Muhanna and D. J. Hallenbeck}, Subordination families and extreme points.
\textit{Trans. Amer. Math. Soc.} \textbf{308} (1988), 83--89.


\bibitem{ACP} J. M. Anderson, J. Clunie and C. Pommerenke, On Bloch
functions and normal functions. \textit{J. Reine Angew. Math.} \textbf{270} (1974), 12--37.

\bibitem{Bo} M. Bonk, The support points of the unit ball in Bloch
space. \textit{J. Funct. Analysis} \textbf{123} (1994), 318--335.

\bibitem{BCK} D. A. Brannan, J. G. Clunie and W. E. Kirwan, On the
coefficient problem for functions of bounded boundary rotation.
\textit{Ann. Acad. Sci. Fenn. A. I.} \textbf{523} (1973), 1--18.


\bibitem{BHHW} L. Brickman, D. J. Hallenbeck, T. H. MacGregor and D. R.
Wilken, Convex hulls and extreme points of families of starlike and
convex mappings.
\textit{Trans. Amer. Math. Soc.} \textbf{185} (1973), 413--428.


\bibitem{CW} J. A. Cima and W. R. Wogen, Extreme points of
the unit ball of the Bloch space $\mathscr{B}_0$. \textit{Michigan
Math. J.} \textbf{25} (1978), 213--222.

\bibitem{Co} F. Colonna, The Bloch constant of bounded harmonic
mappings. \textit{ Indiana Univ. Math. J.} \textbf{38} (1989), 829--840.

\bibitem{Co1}F. Colonna, Bloch and normal functions and their
relation. \textit{Rend. Circ. Mat. Palermo} II \textbf{38} (1989),
161--180.

\bibitem{Du2} { P. Duren}, Univalent Functions. \textit{Springer-Verlag}, New York, Berlin, Heidelberg, Tokyo, 1982.

\bibitem{Pe}{ P. Duren}, Harmonic mappings in the plane. \textit{Cambridge university Press}, New York, 2004.

\bibitem{HH} { D. J. Hallenbeck and K. T. Hallenbeck}, Classes of
analytic functions subordinate to convex functions and extreme
points. \textit{J. Math. Anal. Appl.} \textbf{282} (2003), 792--800.

\bibitem{HM} { D. J. Hallenbeck and T. H. Macgregor},
Subordination and extreme-point theory. \textit{Pacific J. Math.}
\textbf{50} (1974), 455--468.

\bibitem{KW-2019} {I. R. Kayumov and K.-J. Wirths}, 
On the sum of squares of the coefficients of Bloch functions.
\textit{Monatsh Math.} \textbf{190}(1) (2019), 123--135.
 
\bibitem{KW-2019b} {I. R. Kayumov and K.-J. Wirths}, 
Coefficients problems for Bloch functions.
\textit{Anal.Math.Phys.} (2019). https://doi.org/10.1007/s13324-019-00303-z

\bibitem{Mac}T. H. MacGregor, Applications of extreme point theory
to univalent functions. \textit{Michigan Math. J.} \textbf{19} (1972),
361--376.

\bibitem{Mac1}T. H. MacGregor, Hull subordination and extremal
problems for starlike and spirallike mappings. \textit{Trans. Amer.
Math. Soc.} \textbf{183} (1973), 499--510.

\bibitem{RW} St. Ruscheweyh, and  K.-J. Wirths,
On extreme Bloch functions with prescribed critical points. 
\textit{Math. Z.} \textbf{180}(1982), 91--105.

\bibitem{SS} T. Shell-Small, Complex Polynomials. \textit{Cambridge University Press}, New York, 2002.

\bibitem{ST} T. Sugawa and T. Terada, 
A coefficient inequality for Bloch functions with application to uniformly locally univalent functions. 
\textit{Monatshefte Math.} \textbf{156} (2009), 167--173.


\bibitem{Tk}{ K. Tkaczy\'{n}ska}, On extreme points of subordination
families with a convex majorant. \textit{J. Math. Anal. Appl.} \textbf{145} (1990), 216--231.

\bibitem{WirthXiao}  K.-J. Wirths  and J. Xiao, Recognizing $Q_{p,0}$   
functions per Dirichlet space structure.
\textit{Bull. Belg. Math. Soc. Simon Stevin} \textbf{8}(1) (2001), 47--59.


\bibitem{ZS}O. Zariski and P. Samuel, Commutative Algebra, Vol. II.
The University Series in Higher Mathematics. \textit{D. Van Nostrand
Co.}, Inc., Priceton, N. J.-Toronto-London-New York, 1960.

\end{thebibliography}
\end{document}